\documentclass[10pt,a4paper]{article}

\usepackage[T1]{fontenc}
\usepackage[latin1]{inputenc}
\usepackage{graphics}
\usepackage[dvips]{graphicx}
\usepackage[dvips]{color}
\usepackage{amssymb}
\usepackage{amsmath}
\usepackage{amsthm}
\usepackage{epsfig}
\usepackage{fancyhdr}
\usepackage{t1enc,a4,array}
\usepackage[finnish,english]{babel}
\usepackage{mathrsfs}
\usepackage{mathabx}

\defineshorthand{ä}{\"a}
\defineshorthand{Ä}{\"A}
\defineshorthand{ö}{\"o}
\defineshorthand{Ö}{\"O}

\title{}
\author{Heikki Tikanmäki}

\begin{document}
\linespread{1.1}

\newtheorem{maar}{Definition}[section]
\newtheorem{lemma}[maar]{Lemma}
\newtheorem{huom}[maar]{Remark}
\newtheorem{vasta}[maar]{Counter Example}
\newtheorem{theor}[maar]{Theorem}
\newtheorem{corollary}[maar]{Corollary}
\newtheorem{example}[maar]{example}
\newtheorem{prop}[maar]{Proposition}

\numberwithin{equation}{section}

\newcommand{\E}{\mathbb{E}}   
\newcommand{\Pro}{\mathbb{P}} 

\newcommand{\N}{\mathbb{N}}   
\newcommand{\Q}{\mathbb{Q}}   
\newcommand{\R}{\mathbb{R}}   
\newcommand{\Z}{\mathbb{Z}}   
\newcommand{\supp}{\text{supp}}
\newcommand{\Cov}{\text{Cov}}

\begin{center}

\vspace{7 cm}
\Huge

Integral representations of some functionals of fractional Brownian motion

\end{center}

\thispagestyle{empty}

\normalsize

Heikki Tikanmäki, Aalto University, School of Science, P.O. Box 11100, FI-00076 Aalto, heikki.tikanmaki@gmail.com

\vspace{0,5 cm}

30.08.2011

\begin{abstract}
We prove change of variables formulas [It\^o formulas] for functions of both arithmetic and geometric averages of geometric fractional Brownian motion. They are valid for all convex functions, not only for smooth ones. These change of variables formulas provide us integral representations of functions of average in the sense of generalized Lebesgue-Stieltjes integral.



Subject classification (MSC2010): 60G22.
\end{abstract}

\section{Introduction}


In the case of fractional Brownian motion (fBm) it is non-trivial, which functionals have integral representation. In~\cite{a-m-v}, the authors prove that a convex function of the end value of fBm or geometric fBm (gfBm) have an integral representation. In this paper, we prove that analogous integral representations can be constructed also for functionals that are convex functions of the average of fBm or gfBm.

It turns out that the integral representation is the same as for continuous functions of bounded variation. This is not obvious, since for some functionals the integral representations in fBm case and bounded variation case are not the same, see~\cite{a-t-v}.

The usual change of variables formula is not enough for the proofs of the main results, but we need the functional change of variables formula~(\cite{cont}). For the proof of the existence, we use fractional Besov space techniques.

What is a bit surprising here is that we are able to find explicit integral representations for functionals that are functions of the arithmetic average of gfBm. Such formulas are not known even in the case of geometric Brownian motion. The geometric Brownian motion case would correspond to finding Black-Scholes hedging strategies of arithmetic Asian options. 

In this article we will always work in a fixed filtered probability space $(\Omega,\mathcal{F},(\mathcal{F}_t)_{t\geq 0},\Pro)$. The probability space is assumed to be complete. We deal with fractional Brownian motion (fBm) $(B^H(t))_{t\geq 0}$, where the Hurst parameter satisfies $H\in\left(\frac{1}{2},1\right)$. Fractional Brownian motion is a Gaussian process satisfying 
\begin{equation*}
\E B^H(t)=B^H(0)
\end{equation*}
 and
\begin{equation*}
\Cov( B^H(t), B^H(s))=\frac{1}{2}\left ( t^{2H}+s^{2H}-|t-s|^{2H}\right).
\end{equation*}
For this range of $H$, fBm has long-range dependence property. We will mainly consider functionals of geometric fractional Brownian motion (gfBm) $(S(t))_{t\geq0}=\left(\exp B^H(t)\right)_{t\geq 0}$. The stochastic integrals of the article are always pathwise. They are understood in generalized Lebesgue-Stieltjes sense if not mentioned otherwise.

\section{Main results}
\label{results}
\subsection{Existence of stochastic integral}
The aim of this subsection is to prove that the integrands considered in the paper are integrable with respect to fractional Brownian motion in generalized Lebesgue-Stieltjes sense. We proceed here analogously to~\cite{a-m-v}. However, we will need a functional version of the It\^o formula that can be found from~\cite{cont}.

In the proof of the following theorem we will need the concept of fractional Besov spaces. For a short introduction to that topic, see~appendix~A.

Let $f$ be a convex function on $\R$. We refer to~\cite[appendix 3]{revuzyor} for the following facts: The left derivative $f'_-(x)$ (resp. right derivative $f'_+(x)$) exists for all $x\in \R$. If $f$ is differentiable at $x$ then $f'(x)=f'_-(x)$. The second derivative $f''$ of a convex function exists in the sense of distributions and it is a Radon measure on real line. On the other hand, for any Radon measure $\nu$ on $\R$ there is a convex function $g$ such that $g''=\nu$.

In what follows, we will use the following notation. For $t\in[0,T]$ let
\begin{equation*}
G(t)=\exp \left(\frac{1}{T}\int_0^t \log S(s)ds\right)S(t)^{\frac{T-t}{T}},
\end{equation*}
where $S(t)=e^{B^H(t)}$.
\begin{theor}
\label{existence}
Let $f$ be a convex function. Then for $t\in[0,T]$ the integral
\begin{equation*}
\int_0^t \frac{T-s}{T}f'_-(G(s))G(s)dB^H(s)
\end{equation*}
exists almost surely as a generalized Lebesgue-Stieltjes integral.
\end{theor}
The proof is in section~\ref{proofs}. A result similar to theorem~\ref{existence} can be proved analogously for the case of arithmetic average. In that case the theorem takes the following form.
\begin{theor}
\label{existence2}
Let $f$ be a convex function and $t\in[0,T]$, then the integral
\begin{equation*}
\int_0^t f'_-\left (\frac{T-s}{T}S(s)+\frac{1}{T}\int_0^s S(u)du\right )\frac{T-s}{T}S(s)dB^H(s) 
\end{equation*}
exists almost surely as a generalized Lebesgue-Stieltjes integral.
\end{theor}
The proof is in section~\ref{proofs}.
\subsection{Change of variables formulas}
In this subsection the functional change of variables formula of~\cite{cont} is extended for non-smooth convex functions composed with functionals of path of geometric fractional Brownian motion. Considered functionals are geometric and arithmetic averages. For the concepts of horizontal and vertical derivatives and the notation used, see appendix~B.

We begin by proving a change of variables formula (It\^o formula) for the two averages. Note that $dS(t)=S(t)dB^H(t)$.
\begin{prop}
For all $t\in[0,T]$ it holds almost surely that
\label{hedge1}
\begin{align*}
&G(t)=S(0)+\int_0^t\frac{T-s}{T}G(s)dB^H(s),
\end{align*}
where the stochastic integral is understood as a limit of Riemann sums over a sequence of partitions such that the maximum step size goes to zero.
\end{prop}
\begin{corollary}
In particular
\begin{align*}
&\exp{\left (\frac{1}{T}\int_0^T B^H(s)ds\right)}\\
=&S(0)+\int_0^T\frac{T-s}{T}\exp{\left ( \frac{1}{T}\int_0^s B^H(u)du+\frac{T-s}{T}B^H(s)\right )}dB^H(s)
\end{align*}
\end{corollary}
\begin{proof}
Set 
\begin{equation*}
F_t(B^H_t)=S(t)^\frac{T-t}{T}\exp\left(\frac{1}{T}\int_0^tB^H(s)ds\right).
\end{equation*}
The horizontal defivative at time $s\in[0,t]$ of $F$ is
\begin{equation*}
\mathcal{D}_sF_s(B^H_s)=0.
\end{equation*}
The vertical derivative is
\begin{align*}
&\partial_x F(B^H_s)
=\frac{T-s}{T}e^{\frac{1}{T}\int_0^s B^H(u) du+\frac{T-s}{T}B^H(s)}. 
\end{align*}
The second vertical derivative is given analogously by
\begin{equation*}
\partial_x^2 F(B^H_s)=\left (\frac{T-s}{T}\right )^2e^{\frac{1}{T}\int_0^s B^H(u) du+\frac{T-s}{T}B^H(s)}.
\end{equation*}
We know that fractional Brownian motion has zero quadratic variation property for $H>\frac{1}{2}$. Therefore we have by \cite[theorem 3]{cont} that
\begin{align*}
&S(t)^{\frac{T-t}{T}}\exp{\left (\frac{1}{T}\int_0^t B^H(s)ds\right)}\\=&e^{B^H(0)}+\int_0^t\frac{T-s}{T}\exp{\left ( \frac{1}{T}\int_0^s B^H(u)du+\frac{T-s}{T}B^H(s)\right )}dB^H(s),
\end{align*}
where the stochastic integral is understood as a limit of Riemann sums over a sequence of partitions such that the maximum step size goes to zero.
\end{proof}
We have an analogous result for functionals depending on the arithmetic average.
\begin{prop}
\label{arithmavg}
For all $t\in[0,T]$ it holds almost surely that
\begin{equation*}
\frac{T-t}{T}S(t)+\frac{1}{T}\int_0^t S(s)ds=S(0)+\int_0^t \frac{T-s}{T} S(s) dB^H(s),
\end{equation*}
where the stochastic integral is understood as a limit of Riemann sums over a sequence of partitions such that the maximum step size goes to zero.
\end{prop}
\begin{corollary}
In particular
\begin{equation*}
\frac{1}{T}\int_0^T S(s)ds=S(0)+\int_0^T \frac{T-s}{T} S(s) dB^H(s).
\end{equation*}
\end{corollary}
\begin{proof}
The proof goes analogously to the proof of proposition~\ref{hedge1}. Let us define a non-anticipative functional
\begin{equation*}
F_t(B^H_t)=\frac{T-t}{T}S(t)+\frac{1}{T}\int_0^tS(s)ds.
\end{equation*}
The horizontal derivative vanishes and the vertical derivative is given as
\begin{align*}
&\partial_x F(B^H_s)
=\frac{T-s}{T}e^{B^H(s)}.
\end{align*}
The second vertical derivative is also $\frac{T-s}{T}e^{B^H(s)}$. Hence, the change of variables formula of~\cite{cont} takes the form as claimed.
\end{proof}
\begin{huom}
Note that proposition~\ref{arithmavg} could be proved alternatively using integration by parts.
\end{huom}
Now we are ready to provide integral representations first for functionals depending on geometric average and then also for options depending on arithmetic average of gfBm. Finally we obtain corresponding results for arithmetic average when gfBm is replaced by fBm itself.
\begin{theor}
\label{limit1}
Let $t\in[0,T]$, $S(t)=e^{B^H(t)}$ be a geometric fractional Brownian motion with $H\in \left (\frac{1}{2},1\right )$ and $f$ be a convex function. Then it holds almost surely that
\begin{align*}
&f\left(\exp{\left(\frac{1}{T}\int_0^t B^H(s)ds\right)}S(t)^{\frac{T-t}{T}}\right)=&f(S(0))+\int_0^t\frac{T-s}{T} f_-'\left(G(s)\right)G(s)dB^H(s),
\end{align*}
where the stochastic integral in the right side is understood in the sense of generalized Lebesgue-Stieltjes integral.
\end{theor}
The proof is in section~\ref{proofs}.
\begin{corollary}
In particular,
\begin{align*}
&f\left(\exp{\left(\frac{1}{T}\int_0^T B^H(s)ds\right)}\right)
=f(S(0))+\int_0^T \frac{T-s}{T}f_-'\left(G(s)\right)G(s)dB^H(s).
\end{align*}
\end{corollary}
Following theorem is one of the main results, providing integral representations for functionals depending on arithmetic average of gfBm.
\begin{theor}
\label{limit2}
Let $t\in[0,T]$, $S(t)=e^{B^H(t)}$ be a geometric fractional Brownian motion with $H\in \left (\frac{1}{2},1\right )$ and $f$ be a convex function. Then it holds almost surely that
\begin{align*}
&f\left(\frac{T-t}{T}S(t)+\frac{1}{T}\int_0^t S(s)ds\right)\\=&f(S(0))+\int_0^tf'_-\left(\frac{T-s}{T}S(s)+\frac{1}{T}\int_0^sS(u)du\right)\frac{T-s}{T}S(s)dB^H(s),
\end{align*}
where the stochastic integral in the right side is understood in the sense of generalized Lebesgue-Stieltjes integral.
\end{theor}
\begin{proof}
The proof is analogous to the proof of theorem~\ref{limit1}.
\end{proof}
\begin{corollary}
In particular,
\begin{align*}
&f\left(\frac{1}{T}\int_0^T S(s)ds\right)\\=&f(S(0))+\int_0^Tf'_-\left(\frac{T-s}{T}S(s)+\frac{1}{T}\int_0^sS(u)du\right)\frac{T-s}{T}S(s)dB^H(s).
\end{align*}
\end{corollary}
\begin{huom}
The result of theorem~\ref{limit2} can be written also when the geometric fractional Brownian motion $S$ is replaced by a fractional Brownian motion $B^H$ with $H\in\left(\frac{1}{2},1\right)$. In that case we obtain for $t\in [0,T]$ and for a convex function $f$ that
\begin{align*}
&f\left(\frac{T-t}{T}B^H(t)+\frac{1}{T}\int_0^tB^H(s)ds\right)\\=&f(B^H(0))+\int_0^t\frac{T-s}{T}f'_-\left(\frac{T-s}{T}B^H(s)+\frac{1}{T}\int_0^sB^H(u)du\right)dB^H(s)
\end{align*}
almost surely as a generalized Lebesgue-Stieltjes integral.
\end{huom}

\section{Proofs}
\label{proofs}
\subsection{Lemmas}
It is easy to prove the following.
\begin{lemma}
\label{estimatelemma}
Let $\delta \in (0,H)$. Almost surely $\left(S(t)^{1-\frac{t}{T}}\right)_{t\in[0,T]}$ has H\"older continuous sample paths of order $H-\delta$.
\end{lemma}
\begin{lemma}
\label{ubound}
Assume that $X(t)$ is a stochastic process having density $p_t(x)$ s.t. there exists $g(t)\in L^1([0,T])$ s.t. for almost all $t\in(0,T]$ it holds that
\begin{equation*}
p_t(y)\leq g(t)
\end{equation*}
for all $y\in\R$. Let $\alpha\in(0,1)$, then
\begin{equation*}
\E \int_0^T |X(t)+x|^{-\alpha}dt<C<\infty,
\end{equation*}
where $C$ does not depend on $x$.
\end{lemma}
\begin{proof}
First we note that
\begin{equation*}
\E \int_0^T 1_{\{|X(t)+x|\geq 1\}}|X(t)+x|^{-\alpha}dt\leq T.
\end{equation*}
Thus it is enough to consider finiteness of
\begin{align*}
\E \int_0^T 1_{\{|X(t)+x|<1\}}|X(t)+x|^{-\alpha}dt. 
\end{align*}
We note that
\begin{align*}
&\E 1_{\{|X(t)+x|<1\}}|X(t)+x|^{-\alpha}\leq g(t)\int_\R1_{\{|y+x|<1\}}|y+x|^{-\alpha}dy\\\leq& g(t)\int_{x-1}^{x+1}|y+x|^{-\alpha}dy=g(t)\frac{2}{1-\alpha}.
\end{align*}
Now by Fubini's theorem
\begin{align*}
&\E \int_0^T 1_{\{|X(t)+x|<1\}}|X(t)+x|^{-\alpha}dt=\int_0^T \E \left(1_{\{|X(t)+x|<1\}}|X(t)+x|^{-\alpha}\right)dt\\ \leq & \frac{2}{1-\alpha}\|g\|_{L^1([0,T])}.
\end{align*}
\end{proof}
The proof of the following lemma uses Malliavin calculus. For a detailed expression of the topic, see~\cite{nualart}. Let $\mathcal{H}$ be an isonormal Gaussian Hilbert space associated with process $B^H$. Thus $\mathcal{H}$ is equipped with inner product defined by 
\begin{equation*}
<1_{[0,t]},1_{[0,s]}>_\mathcal{H}=\frac{1}{2}(t^{2H}+s^{2H}-|t-s|^{2H})
\end{equation*}
 for $t,s\in[0,T]$. Malliavin derivative of a random variable $F(h)$, $h\in\mathcal{H}$ is denoted by $DF$ and it takes its values in $\mathcal{H}$. The second Malliavin derivative $D^2F$ is an element of space $\mathcal{H}\otimes \mathcal{H}$. Besides of the norms of $\mathcal{H}$ and $\mathcal{H}\otimes \mathcal{H}$ and proper $L^p(\Omega)$ norms, we will need the following norm
\begin{equation*}
\|D^2F\|_{L^p(\Omega; \mathcal{H} \otimes \mathcal{H})}^p=\E (\|D^2F\|^p_{\mathcal{H}\otimes \mathcal{H}}).
\end{equation*}

\begin{lemma}
\label{density1}
Let $X(t)=\frac{T-t}{T}e^{B^H(t)}+\frac{1}{T}\int_0^t e^{B^H(u)}du$. The density of random variable $X(t)$ exists and it is denoted by $p_t(x)$. Furthermore, there exists $g\in L^1([0,T])$ such that $p_t(x)\leq g(t)$ for all $x\in \R$ and almost all $t\in [0,T]$.
\end{lemma}
\begin{proof}
The density of $X(t)$ exists by~\cite[prop.~2.1.1.]{nualart}. By \cite[prop.~2.1.2.]{nualart} we have for $\alpha,\beta, q>0$ s.t. $\frac{1}{\alpha}+\frac{1}{\beta}+\frac{1}{q}= 1$ that
\begin{equation*}
p_t(x)\leq C_{\alpha,\beta,q}\left(\Pro(X(t)>x)^{\frac{1}{q}}\right)\left(\E (\|DX(t)\|_\mathcal{H}^{-1})+\|D^2X(t)\|_{L^\alpha(\Omega;\mathcal{H}\otimes \mathcal{H})}\left\|\|DX(t)\|_\mathcal{H}^{-2}\right\|_{L^\beta(\Omega)}\right).
\end{equation*}
We have that
\begin{equation*}
DX(t)=\frac{T-t}{T}e^{B^H(t)}1_{[0,t]}+\frac{1}{T}\int_0^te^{B^H(u)}1_{[0,u]}du
\end{equation*}
and
\begin{equation*}
D^2 X(t)=\frac{T-t}{T}e^{B^H(t)}1_{[0,t]}\otimes 1_{[0,t]}+\frac{1}{T}\int_0^t e^{B^H(u)}1_{[0,u]}\otimes 1_{[0,u]}du.
\end{equation*}
It holds that
\begin{align}
\label{malliavinnorm}
\|DX(t)\|_\mathcal{H}^2=&\left(\frac{T-t}{T}\right)^2e^{2B^H(t)}t^{2H}\\&+\frac{1}{T^2}\int_0^t\int_0^t e^{B^H(u)+B^H(v)}\frac{1}{2}\left(u^{2H}+v^{2H}-|u-v|^{2H}\right)dudv \nonumber\\&+\frac{T-t}{T^2}\int_0^te^{B^H(t)+B^H(u)}\left(u^{2H}+t^{2H}-(t-u)^{2H}\right)du.\nonumber
\end{align}
Now we have
\begin{align*}
&\int_0^t\int_0^t e^{B^H(u)+B^H(v)}\frac{1}{2}\left(u^{2H}+v^{2H}-|u-v|^{2H}\right)dudv \\\geq& \frac{1}{2}\exp\left({2\inf_{s\in [0,T]}B^H(s)}\right)\int_0^t\int_0^t \left(u^{2H}+v^{2H}-|u-v|^{2H}\right)dudv\\
=&\frac{1}{2}\exp\left({2\inf_{s\in [0,T]}B^H(s)}\right)\left(\frac{2}{2H+1}-\frac{1}{(2H+1)(H+1)}\right)t^{2H+2}.
\end{align*}
Note that all the three terms in the right side of equation~(\ref{malliavinnorm}) are positive. Thus we have for some constant $C_1>0$ that
\begin{equation*}
\|DX(t)\|_\mathcal{H}\geq C_1 \exp{\left(\inf_{s\in [0,T]}B^H(s)\right)}\left(\frac{T-t}{T}t^H\vee t^{H+1}\right).
\end{equation*}
The norm of the second Malliavin derivative can be bounded as
\begin{align*}
&\|D^2 X(t)\|_{\mathcal{H}\otimes \mathcal{H}}^2=\left(\frac{T-t}{T}\right)^2 e^{2B^H(t)}<1_{[0,t]}\otimes 1_{[0,t]},1_{[0,t]}\otimes 1_{[0,t]}>_{\mathcal{H}\otimes \mathcal{H}}\\
&+2\frac{T-t}{T^2}\int_0^te^{B^H(u)+B^H(t)}<1_{[0,t]}\otimes 1_{[0,t]},1_{[0,u]}\otimes 1_{[0,u]}>_{\mathcal{H}\otimes \mathcal{H}}du\\
&+\frac{1}{T^2}\int_0^t\int_0^te^{B^H(u)+B^H(v)}<1_{[0,u]}\otimes 1_{[0,u]},1_{[0,v]}\otimes 1_{[0,v]}>_{\mathcal{H}\otimes \mathcal{H}}dudv\\
=&\left(\frac{T-t}{T}\right)^2 e^{2B^H(t)}t^{4H}+2\frac{T-t}{T^2}\int_0^t e^{B^H(u)+B^H(t)}\frac{1}{4}\left(t^{2H}+u^{2H}-(t-u)^{2H}\right)^2du\\
&+\frac{1}{T^2}\int_0^t\int_0^t e^{B^H(u)+B^H(v)}\frac{1}{4}\left(v^{2H}+u^{2H}-|v-u|^{2H}\right)^2dudv\\
\leq &\tilde{C}_2\exp{\left(2\sup_{s\in [0,T]}B^H(s)\right)}\left(t^{4H}+t^{4H+2}\right)\leq C_2\exp{\left(2\sup_{s\in [0,T]}B^H(s)\right)}t^{4H},
\end{align*}
for some constants $\tilde{C}_2$ and $C_2$. Now,
\begin{align*}
p_t(x)\leq&C_{\alpha,\beta,q}
(C_1^{-1}\E \exp{\left(-\inf_{s\in [0,T]}B^H(s)\right)}\left(t^{-H}\frac{T}{T-t}\wedge t^{-H-1}\right) \\
&+C_1^{-2\beta}C_2\E e^{\alpha \sup_{s\in[0,T]}B^H(s)}\E e^{-\beta \inf_{s\in [0,T]}B^H(s)}t^{4H\frac{\alpha}{2}}t^{-{(2H+2)}\frac{\beta}{2}})\\
\leq & C\left( \left(t^{-H}\frac{T}{T-t}\wedge t^{-H-1}\right)+t^{2H\alpha}t^{-{(H+1)}\beta}\right)\in L^1([0,T])
\end{align*}
if we choose $\alpha,\beta,q$ such that $2\alpha H-(H+1)\beta>-1$. For example $\alpha=q=4$ and $\beta=2$ is a possible choice.
\end{proof}
\begin{lemma}
\label{density2}
The density of random variable $\log G(t)$ exists and is denoted by $p_t(x)$. Furthermore, there exists $g \in L^1([0,T])$ such that $p_t(x)\leq g(t)$ for all $x\in \R$ and almost all $t\in [0,T]$.
\end{lemma}
\begin{proof}
The existence of $p_t(x)$ is due to Gaussianity of $\log G(t)$.
\begin{equation*}
\log G(t)=\frac{T-t}{T}B^H(t)+\frac{1}{T}\int_0^t B^H(u)du.
\end{equation*}
The Malliavin derivative is
\begin{equation*}
D\log G(t)=\frac{T-t}{T}1_{[0,t]}+\frac{1}{T}\int_0^t 1_{[0,u]}du.
\end{equation*}
The second Malliavin derivative $D^2 \log G(t)=0$. Hence, by \cite[prop.~2.1.2.]{nualart}
\begin{equation*}
p_t(x)\leq C \left(t^{-H}\frac{T}{T-t}\wedge t^{-H-1}\right) \in L^1([0,T]).
\end{equation*}
\end{proof}
\begin{lemma}
\label{yes1}
Let $p\geq 1$. Then there exists $C<\infty$ such that
\begin{equation*}
\E |\log G(t)-\log G(s)|^p\leq C|t-s|^{pH}.
\end{equation*}
\end{lemma}
\begin{proof}
Follows from the corresponding property of fBm.
\end{proof}
\begin{lemma}
Let $p\geq 1$ and
\begin{equation*}
X(t)=\frac{T-t}{T}e^{B^H(t)}+\frac{1}{T}\int_0^te^{B^H(u)}du.
\end{equation*}
 Then there exists $C<\infty$ such that
\begin{equation*}
\E \left | X(t)-X(s)\right|^p\leq C|t-s|^{pH}.
\end{equation*}
\label{yes2}
\end{lemma}
\begin{proof}
\begin{align*}
&|X(t)-X(s)|\leq
\left| e^{B^H(t)}-e^{B^H(s)}\right|\\&+\frac{1}{T}\left| t e^{B^H(t)}-te^{B^H(s)}\right|+\frac{1}{T}\left|te^{B^H(s)}-se^{B^H(s)}\right|+\frac{|t-s|}{T}\sup_{u\in[0,T]}e^{B^H(u)}
\\\leq& 2\left|e^{B^H(t)}-e^{B^H(s)}\right|+|t-s|\frac{1}{T}\left(e^{B^H(s)}+\sup_{u\in [0,T]}e^{B^H(u)} \right).
\end{align*}
Note that for $a,b\geq 0$ it holds that $(a+b)^p\leq 2^p(a^p+b^p)$. Thus it is enough to show the claim term-wise. The last term is obvious. For the first one we note that
\begin{align*}
|e^x-e^y|=\left|\int_x^y e^udu\right|\leq |x-y|e^{x\vee y}.
\end{align*}
Now we have by H\"older inequality that
\begin{align*}
&\E \left | e^{B^H(t)}-e^{B^H(s)}\right |^p\\\leq& \sqrt{\E |B^H(t)-B^H(s)|^{2p}}\sqrt{\E \sup_{u\in [0,T]}|e^{B^H(u)}|^{2p}}\\
\leq & \sqrt{\tilde{C}|t-s|^{2pH}}\sqrt{\E e^{2p\sup_{u\in [0,T]} B^H(u)}}\leq C |t-s|^{pH}.
\end{align*}
Finiteness of ${\E e^{2p\sup_{u\in [0,T]} B^H(u)}}$ follows from \cite[p. 182-184]{lifshits}.

\end{proof}
\subsection{Proofs of main results}
\begin{proof}[Proof of theorem~\ref{existence}]
We have to show that for some $\beta \in (1-H,\frac{1}{2})$ it holds that
\begin{equation*}
\left|\left|\frac{T-s}{T}f'_-(G(s))G(s)\right|\right|_{2,\beta}<\infty \quad \text{almost surely}.
\end{equation*}
Then we obtain the claim using theorem~\ref{integralexists}. First of all,
\begin{equation*}
\int_0^T \left | \frac{T-s}{T}f'_-(G(s))G(s)\right |\frac{1}{s^\beta}ds\leq \sup_{u\in[0,T]} \left |f'_-(G(u))G(u)\right|\int_0^T s^{-\beta}ds<\infty \quad a.s..
\end{equation*}
The other term of the Besov norm is more complicated. Define
\begin{equation*}
g(t)=\frac{T-t}{T}f'_-(G(t))G(t).
\end{equation*}
We obtain for $0\leq s \leq t \leq T$ that
\begin{align*}
&|g(t)-g(s)|=
\frac{1}{T}\exp{\left ( \frac{1}{T}\int_0^s \log{S(u)}du\right )}\\ &\times\left | \exp{\left ( \frac{1}{T}\int_s^t \log{S(u)}du\right )}(T-t)S(t)^{1-\frac{t}{T}}f'_-(G(t))-(T-s)S(s)^{1-\frac{s}{T}}f'_-(G(s))\right |.
\end{align*}
For the difference we obtain by the triangle inequality that
\begin{align*}
&\left| \exp{\left ( \frac{1}{T}\int_s^t \log{S(u)}du\right )}(T-t)S(t)^{1-\frac{t}{T}}f'_-(G(t))-(T-s)S(s)^{1-\frac{s}{T}}f'_-(G(s))\right |\\
\leq & \left|\exp{\left ( \frac{1}{T}\int_s^t \log{S(u)}du\right )}-1\right|(T-t)S(t)^{1-\frac{t}{T}}\left|f'_-(G(t))\right|\\
&+(t-s)S(t)^{1-\frac{t}{T}}|f'_-(G(t))|\\
&+(T-s)\left |S(t)^{1-\frac{t}{T}}-S(s)^{1-\frac{s}{T}}\right |\left|f'_-(G(t))\right|\\
&+(T-s)S(s)^{1-\frac{s}{T}}\left | f'_-(G(t))-f'_-(G(s))\right |\\
=:&A_1+A_2+A_3+A_4.
\end{align*}
Let us proceed term-wise:
\begin{align*}
&\left|\exp{\left ( \frac{1}{T}\int_s^t \log{S(u)}du\right )}-1\right|\leq \sum_{k=1}^\infty \frac{\left | \frac{1}{T}\int_s^t \log S(u)du\right |^k}{k!}\\
\leq & \frac{1}{T}(t-s)\sup_{u\in[0,T]} |\log S(u)|\sum_{k=0}^\infty \frac{\left |\frac{1}{T}(t-s)\sup_{u\in[0,T]} |\log S(u)| \right |^k}{(k+1)!}\\
\leq & \frac{1}{T}(t-s)\sup_{u\in[0,T]} |\log S(u)| \exp{\left ( \frac{1}{T}(t-s)\sup_{u\in[0,T]}|\log S(u)|\right )}.
\end{align*}
Hence,
\begin{align*}
&\int_0^T\int_0^t \frac{A_1}{(t-s)^{\beta+1}}dsdt\\\leq&\sup_{u\in[0,T]}S(u)^{1-\frac{u}{T}}\sup_{u\in[0,T]}|f'_-(G(u))|\sup_{u\in[0,T]}|\log S(u)|e^{\sup_{u\in[0,T]}|\log S(u)|}\\&\times\int_0^T \int_0^t (t-s)^{-\beta}dsdt<\infty
\end{align*}
almost surely. For the second term we have
\begin{align*}
&\int_0^T \int_0^t \frac{A_2}{(t-s)^{\beta+1}}dsdt\\\leq& \left (\sup_{u\in[0,T]}S(u)^{1-\frac{u}{T}}\right )\sup_{u\in[0,T]}|f'_-(G(u))|\int_0^T\int_0^t (t-s)^{-\beta}dsdt<\infty
\end{align*}
almost surely. For the third one we can use lemma~\ref{estimatelemma}. Let us choose $\delta=\frac{H-\beta}{2}$. Now $0<\delta<H-\beta$. Thus there exists almost surely finite constant $C(\omega)$ such that
\begin{align*}
&\int_0^T \int_0^t \frac{A_3}{(t-s)^{\beta+1}}dsdt \leq T \sup_{u\in[0,T]}|f'_-(G(u))|C(\omega) \int_0^T \int_0^t(t-s)^{(H-\beta)-\delta-1}dsdt<\infty
\end{align*}
almost surely. We need an estimate for the fourth term and then we are done. Let us denote the second (distribution) derivative of $f$ by $\mu$. Let us assume first that $\kappa=\supp(\mu)$ is compact. By \cite[p. 545]{revuzyor} or \cite{a-m-v}, we have estimate
\begin{align*}
&\int_0^T\int_0^t \frac{A_4}{(t-s)^{\beta+1}}dsdt\\
\leq& T \left (\sup_{u\in [0,T]}S(u)^{1-\frac{u}{T}}\right)\int_0^T \int_0^t\int_\kappa\frac{1_{\{G(s)<x<G(t)\}}+1_{\{G(t)<x<G(s)\}}}{(t-s)^{\beta+1}}\mu(dx)dsdt.
\end{align*}
Consider now integral of the first indicator. The other one can be considered analogously.
\begin{equation*}
J=\int_0^T \int_0^t \int_\kappa\frac{1_{\{G(s)<x<G(t)\}}}{(t-s)^{\beta+1}}\mu(dx)dsdt.
\end{equation*}
By Tonelli's theorem we have that
\begin{align*}
\E J=\int_\kappa \E \left( \int_0^T \int_0^t \frac{1_{\{G(s)<x<G(t)\}}}{(t-s)^{\beta+1}}dsdt\right)\mu(dx).
\end{align*}
Let us define now
\begin{align*}
T_t(x):=&\sup \left \{u\in [0,t]:G(u)=x\right\},
\end{align*}
with the convention that supremum over an empty set is $0$. On the set $\{\omega\in \Omega : x<G(t)\}$ it holds that $T_t(x)<t$ a.s.. Thus,
\begin{align*}
&\int_0^t\frac{1_{\{G(s)<x<G(t)\}}}{(t-s)^{\beta+1}}ds\\ \leq& \int_0^{T_t(x)}\frac{1_{\{x<G(t)\}}}{(t-s)^{\beta+1}}ds=1_{\{x<G(t)\}}\frac{(t-T_t(x))^{-\beta}-t^{-\beta}}{\beta}.
\end{align*}
In the case that $T_t(x)=0$, this upperbound is zero. In what follows we assume that $0<T_t(x)<t$. We define process $(Y(t))_{t\in[0,T]}$ by
\begin{equation*}
Y(t)=\log G(t).
\end{equation*}
We use lemma~\ref{yes1} and Garsia-Rodemich-Rumsey inequality~\cite{garsia} or~\cite[A.3]{nualart}, and obtain that for $p\geq 1$ and $\gamma\in(\frac{1}{p},H)$ there exists $D(\omega)$ such that for all $s,t\in[0,T]$
\begin{equation*}
\left | Y(t)-Y(s)\right |^p\leq D(\omega)|t-s|^{\gamma p-1}\int_0^T \int_0^T\frac{\left| Y(u_1)-Y(u_2)\right |^p}{|u_1-u_2|^{\gamma p +1}}du_2du_1.
\end{equation*}
Substituting $s=T_t(x)$ to the inequality we obtain
\begin{align*}
&\left | Y(t)-\log x \right |^p\\\leq& D(\omega)(t-T_t(x))^{\gamma p -1}\int_0^T \int_0^T\frac{\left| Y(u_1)-Y(u_2)\right |^p}{|u_1-u_2|^{\gamma p +1}}du_2du_1.
\end{align*}
Fix $\epsilon \in (0,H-\frac{1}{2})$ and choose $p=\frac{2}{\epsilon}$ and $\gamma=H-\frac{\epsilon}{2}$. Let $\xi \in(0,\frac{\epsilon}{2})$. By $H-\xi$ -H\"older continuity of paths of $Y$ it holds that
\begin{equation}
\label{Hbound}
\int_0^T \int_0^T\frac{\left| Y(u_1)-Y(u_2)\right |^\frac{2}{\epsilon}}{|u_1-u_2|^{\frac{2H}{\epsilon}}}du_2du_1\leq \tilde{C} \int_0^T \int_0^T |u_1-u_2|^{-\frac{2\xi}{\epsilon}}du_1du_2<\infty \quad \text{a.s.}.
\end{equation}
Now in the set $\{\omega \in \Omega : G(s)<x<G(t)\}\subset\{\omega \in \Omega : x<G(t)\}\cap \{\omega\in\Omega : T_t(x)>0\}$ we have that
\begin{align*}
&\left | Y(t) - \log x\right |\\\leq& D(\omega)^\frac{\epsilon}{2}(t-T_t(x))^{H-\epsilon}\left ( \int_0^T \int_0^T\frac{\left| Y(u_1)-Y(u_2)\right |^\frac{2}{\epsilon}}{|u_1-u_2|^{\frac{2H}{\epsilon}}}du_2du_1\right )^\frac{\epsilon}{2}.
\end{align*}
For each $t>0$, it holds almost surely that
\begin{equation*}
\left | Y(t) - \log x \right |>0.
\end{equation*}
Hence,
\begin{align*}
(t-T_t(x))^{-\beta}\leq &D(\omega)^{\frac{\beta \epsilon}{2(H-\epsilon)}}\left|Y(t)-\log x\right|^{-\frac{\beta}{H-\epsilon}}\\&\times\left ( \int_0^T \int_0^T\frac{\left| Y(u_1)-Y(u_2)\right |^\frac{2}{\epsilon}}{|u_1-u_2|^{\frac{2H}{\epsilon}}}du_2du_1\right )^\frac{\epsilon \beta}{2(H-\epsilon)}
\end{align*}
almost surely. Using equation~(\ref{Hbound}), we have for some almost surely finite random variable $D_2(\omega)$ that
\begin{equation*}
\int_0^T (t-T_t(x))^{-\beta}dt\leq D_2(\omega)\int_0^T \left | Y(t) -\log x\right |^{-\frac{\beta}{H-\epsilon}}dt.
\end{equation*}
Note that $D_2(\omega)$ does not depend on $x$. Now we use lemmas~\ref{ubound} and~\ref{density2} 
and obtain that 
\begin{equation*}
\E \left ( \int_0^T \left | Y(t)- \log x\right |^{-\frac{\beta}{H-\epsilon}}dt\right )<c_1<\infty.
\end{equation*}
This implies that
\begin{equation}
\label{c2}
\E \left( \int_0^T \int_0^t \frac{1_{\{G(s)<x<G(t)\}}}{(t-s)^{\beta+1}}dsdt\right)<c_2<\infty,
\end{equation}
where $c_1, c_2$ do not depend on $x$. Thus,
\begin{equation*}
\int_0^T\int_0^t \frac{A_4}{(t-s)^{\beta+1}}dsdt<\infty \quad a.s.,
\end{equation*}
because $\mu(\kappa)<\infty$. This proves the theorem in the case that $\supp(\mu)$ is compact.

Let us now consider the general case with no assumption on the compactness of $\supp(\mu)$. For any $n\in\Z^+$ we define
\begin{equation}
\label{omegapartition}
\Omega_n=\left\{\omega \in \Omega : \left(\max_{u\in [0,T]}G(u)\right) \in [0,n]\right\}.
\end{equation}
Let us define now a new convex function
\begin{equation*}
f_n(x)=f(x)1_{\{x\in[0,n]\}}+(f'_-(n)(x-n)+f(n))1_{\{x>n\}}+(f'_+(0)x+f(0))1_{\{x<0\}}.
\end{equation*}
Now $f_n(x)=f(x)$, when $x\in[0,n]$. Let us denote the second derivative of $f_n$ by $\mu_n$. We know that $\supp(\mu_n)\subset [0,n]$ is compact. We know also that almost surely
\begin{equation*}
\int_0^t f_n'(G(s))\frac{T-s}{T}G(s)dB^H(s)
\end{equation*}
is well defined a.s. on $\Omega_n$. By the definition of sets $\Omega_n$, we know that $\Omega=\bigcup_{n=1}^\infty \Omega_n$. Thus, the stochastic integral is well-defined in generalized Lebesgue-Stieltjes sense almost surely for $\omega \in \Omega$. This completes the proof of the theorem.
\end{proof}
\begin{proof}[Proof of theorem~\ref{existence2}]
The proof is analogous to the proof of theorem~\ref{existence}. We need lemma~\ref{yes2} to apply Garsia-Rodemich-Rumsey theorem to process
\begin{equation*}
X(t)=\frac{T-t}{T}e^{B^H(t)}+\frac{1}{T}\int_0^t e^{B^H(u)}du.
\end{equation*}
After that we use lemma~\ref{density1} instead of lemma~\ref{density2}.
\end{proof}

\begin{proof}[Proof of theorem~\ref{limit1}]
First we prove the result for smooth functions $f$ and then we use fractional Besov space techniques to obtain the general case.

If $f\in C^2$, then we can use the It\^o formula of~\cite{cont} and lemmas~\ref{verticalchain},~\ref{productrule} and~\ref{horizontalchain} to  obtain that
\begin{align}
\label{ito1}
&f\left ( e^{\frac{1}{T}\int_0^t B^H(s)ds}S(t)^\frac{T-t}{T}\right )\\=&f\left (e^{B^H(0)}\right)+\int_0^t f'_-\left(\exp\left(\frac{1}{T}\int_0^s B^H(u)du+\frac{T-s}{T}B^H(s)\right )\right )\frac{T-s}{T}\nonumber\\&\times\exp\left(\frac{1}{T}\int_0^s B^H(u)du+\frac{T-s}{T}B^H(s)\right )dB^H(s),
\nonumber\\
=&f\left(e^{B^H(0)}\right)+\int_0^t \frac{T-s}{T}f'_-(G(s))G(s)dB^H(s).\nonumber
\end{align}
where the stochastic integral is understood as a limit of Riemann sums. The aim from now on is to show that equation~(\ref{ito1}) holds for any convex $f$ and the integral is generalized Lebesgue-Stieltjes integral.

Let now $f$ be convex function, $f''=\mu$ in the sense of distributions and $\phi$ some positive $C^\infty((-\infty,0])$ function with compact support such that
\begin{equation*}
\int_{-\infty}^0 \phi(y)dy=1.
\end{equation*}
Define now for $n\in \Z_+$ the approximating functions
\begin{equation*}
f_n(x)=n\int_{-\infty}^0 f(x+y)\phi(ny)dy.
\end{equation*}
For each $n$, $f_n\in C^\infty$ is convex and locally bounded (\cite[p. 221]{revuzyor}). It also holds for $g\in C^\infty$ with compact support that
\begin{equation}
\label{weakconvergence}
\lim_{n\rightarrow \infty} \int_\R g(x)f_n''(x)dx=\int_\R g(x)\mu(dx).
\end{equation}
Define now for $k\in \N $, $\Omega_k$ as in equation~(\ref{omegapartition}). Let us choose some convex functions $f_{n,k}\in C^2(\R)$ s.t. $f_{n,k}(x)=f_n(x)$, when $x\in [0,k]$ and $f_{n,k}''(x)=0$, when $x\in \R \backslash [-1,k+1]$ s.t. 
\begin{equation*}
\int_{-1}^0 f_{n,k}''(x)dx+\int_{k}^{k+1}f_{n,k}''(x)dx\leq 1.
\end{equation*}
Note that convexity is preserved as long as $f_{n,k}''\geq 0$. Now $f_{n,k}$ is a convex $C^2$ function with compactly supported second derivative. By equation~(\ref{ito1}) we have that the change of variables formula holds for functions $f_n$ and $f_{n,k}$ for all $k,n\in \N$. It holds almost surely by~\cite[p.221]{revuzyor} that
\begin{align*}
&f_{n}\left(G(t))\right)\rightarrow f\left(G(t)\right)
\end{align*}
and
\begin{equation*}
f_{n}'(G(t)) \rightarrow f'_-(G(t)).
\end{equation*}
For the convergence of the stochastic integral we will use the theory of appendix~A. It turns out that it is sufficient to show for some $\beta \in \left(1-H,\frac{1}{2}\right)$ that almost surely
\begin{align*}
\left|\left|f'_n(G(t))G(t)-f'_-(G(t))G(t)\right|\right|_{2,\beta}\rightarrow 0, \quad \text{as}\quad n\rightarrow \infty.
\end{align*}
We will first show this in set $\Omega_k$. Fix $\beta \in\left(1-H,\frac{1}{2}\right)$. We have by convexity that
\begin{equation*}
\sup_{u\in [0,T]}|f'_n(G(u))|=\left|f_n'\left(\sup_{u\in [0,T]}G(u)\right)\right|\vee\left|f_n'\left(\inf_{u \in [0,T]}G(u)\right)\right|.
\end{equation*}
Thus, we have for $n$ large enough that
\begin{align*}
&\frac{\left| f_n'(G(t))G(t)-f'_-(G(t))G(t)\right |}{t^\beta}\\\leq& \frac{2\sup_{u\in[0,T]}\left(G(u)\right)\left(\sup_{u\in[0,T]}\left|f'_-(G(u))\right|+1\right)}{t^\beta}\in L^1([0,T]).
\end{align*}
The Lebesgue dominated convergence theorem now implies that almost surely
\begin{equation*}
\int_0^T \frac{\left | f_n'(G(t))G(t)-f'_-(G(t))G(t)\right |}{t^\beta}dt \rightarrow 0.
\end{equation*}
For the other term of the norm we obtain 
\begin{align*}
&\frac{|f_n'(G(t))G(t)-f'_-(G(t))G(t)-(f_n'(G(s))G(s)-f'_-(G(s))G(s))|}{(t-s)^{\beta+1}}\\
\leq &\frac{|f_n'(G(t))G(t)-f_n'(G(s))G(s)|}{(t-s)^{\beta+1}}+\frac{|f'_-(G(t))G(t)-f_-'(G(s))G(s)|}{(t-s)^{\beta+1}}.
\end{align*}
We have
\begin{align}
\label{bound}
&\frac{|f_n'(G(t))G(t)-f_n'(G(s))G(s)|}{(t-s)^{\beta+1}}\\
\leq & \frac{|f_{n}'(G(t))G(t)-f_{n}'(G(t))G(s)|}{(t-s)^{\beta+1}}+\sup_{u\in[0,T]}G(u)\frac{|f_{n,k}'(G(t))-f_{n,k}'(G(s))|}{(t-s)^{\beta+1}}.
\nonumber
\end{align}
Choose $\delta=\frac{H-\beta}{2}$. We have for $n$ large enough using lemma~\ref{estimatelemma} that
\begin{align*}
&\int_0^T \int_0^t\frac{|f_{n}'(G(t))G(t)-f_{n}'(G(t))G(s)|}{(t-s)^{\beta+1}}dsdt\\\leq & C(\omega)\int_0^T \int_0^t\left(\sup_{u\in[0,T]}\left | f'_-\left( G(u)\right)\right |+1\right)\frac{(t-s)^{H-\delta}}{(t-s)^{\beta+1}}dsdt<\infty
\end{align*}
and by the proof of theorem~\ref{existence} that
\begin{align*}
&\E \int_0^T \int_0^t\frac{|f_{n,k}'(G(t))-f_{n,k}'(G(s))|}{(t-s)^{\beta+1}}dsdt\\\leq& \int_{-1}^{k+1}c_2 f_{n,k}''dx\leq c_2+c_2\int_0^k f_n''(x)dx<C<\infty,
\end{align*}
where $c_2$ is as in equation~(\ref{c2}) and $C$ does not depend on $n$ because $\int_0^k f_n''(x)dx\rightarrow \mu([0,k])$ by equation~(\ref{weakconvergence}). Note that $\mu([0,k])<\infty$, because $\mu$ is a Radon measure.

Now we get by dominated convergence theorem that almost surely in $\Omega_k$
\begin{align*}
&\int_0^T\int_0^t \frac{\left | f'_n(G(t))G(t)-f'_n(G(s))G(s)-\left( f'_-(G(t))G(t)-f'_-(G(s))G(s)\right )\right |}{(t-s)^{\beta+1}}dsdt\\&\rightarrow_{n\rightarrow \infty}0
\end{align*}
and thus almost surely in $\Omega_k$
\begin{equation}
\label{conv}
||f_n'(G(t))G(t)-f'_-(G(t))G(t)||_{2,\beta}\rightarrow_{n\rightarrow \infty}0.
\end{equation}
We note that $\Omega=\bigcup_{k=1}^\infty \Omega_k$ and thus~(\ref{conv}) holds almost surely in $\Omega$. This implies that the approximating integrals converge to an integral in generalized Lebesgue-Stieltjes sense. Moreover, the limiting integral is what is claimed.
\end{proof}
\section{Conclusions}
In this paper we were able to extend the functional It\^o formula of~\cite{cont} for non-smooth convex functions in the special case of driving gfBm or fBm and functional depending on the average of the driving process.





The results of section~\ref{results} remain true if we add such deterministic drift to fBm that does not change path properties. That is, we can add H\"older continuos drift with zero quadratic variation.

For a smooth $f$ the integral representations are limits of Riemann sums. This can be used as a starting point for developing a discretizing method for the stochastic integral. 
\subsection{Financial interpretation}
The results of the paper can be used for obtaining hedges for Asian options in fractional Black-Scholes model. In this model, the stock price is modeled by geometric fBm $S(t)=\exp B^H(t)$. However, this model allows for arbitrage opportunities (\cite{cheridito, shiryaev, a-m-v}). Thus, the use of such model as a financial model is questionable. In fact the results of this paper can be used for obtaining new arbitrage examples in that model.
\subsubsection{Arbitrage opportunity}
Theorems~\ref{limit1}~and~\ref{limit2} provide us concrete examples of arbitrage opportunities in fractional Black-Scholes model. In the setup of Theorem~\ref{limit2} we say that the arithmetic Asian call option with payoff
\begin{equation*}
\left ( \frac{1}{T}\int_0^T S(s)ds-K\right )^+
\end{equation*}
is out-of-the-money at time $t<T$ if
\begin{equation*}
\frac{1}{T}\int_0^tS(s)ds+\frac{T-t}{T}S(t)<K.
\end{equation*}
It is easy to see that out-of-the-money options in this sense have zero hedging cost, even though the option will be in-the-money at maturity $T$ with positive probability.
\section*{Appendix A: Fractional Besov spaces and generalized Lebesgue-Stieltjes integrals}
For the following facts, we refer to~\cite{mishura}. Let $\beta\in(0,1)$.
\begin{maar}
Fractional Besov space of order $1,\beta$ is denoted by $W_1^\beta([0,T])$. Let $f:[0,T]\mapsto \R$ be measurable. We say that $f\in W_1^\beta([0,T])$ if
\begin{equation*}
||f||_{1,\beta}=\sup_{0\leq s< t\leq T}\left(\frac{|f(t)-f(s)|}{(t-s)^\beta}+\int_s^t\frac{|f(u)-f(s)|}{(u-s)^{\beta+1}}du\right)<\infty.
\end{equation*}
\end{maar}
Note that $||\cdot||_{1,\beta}$ is not a norm but a seminorm.
\begin{maar}
Fractional Besov space of order $2,\beta$ is denoted by $W_2^\beta([0,T])$. Let $f:[0,T]\mapsto \R$ be measurable. We say that $f\in W_2^\beta([0,T])$ if
\begin{equation*}
||f||_{2,\beta}=\int_0^T\frac{|f(t)|}{t^\beta}dt+\int_0^T\int_0^t \frac{|f(t)-f(s)|}{(t-s)^{\beta+1}}dsdt<\infty.
\end{equation*}
\end{maar}
Let us denote by $C^\alpha([0,T])$ the H\"older continuous functions of order $\alpha$ on $[0,T]$. It holds for $\epsilon\in(0,\beta\wedge(1-\beta))$ that
\begin{equation*}
C^{\beta+\epsilon}([0,T])\subset W_1^\beta([0,T])\subset C^{\beta-\epsilon}([0,T]).
\end{equation*}
Thus for $0<\beta<H$ the trajectories of fBm $B^H$ belong to $W_1^\beta([0,T])$ almost surely. It also holds that
\begin{equation*}
C^{\beta+\epsilon}([0,T])\subset W_2^\beta([0,T]).
\end{equation*}
Let $\alpha>0$. We define Riemann-Liouville fractional integrals on $[0,t]$ as follows
\begin{equation*}
\left(I^\alpha_{0+}f\right)(x)=\frac{1}{\Gamma(\alpha)}\int_0^xf(s)(x-s)^{\alpha-1}ds
\end{equation*}
and
\begin{equation*}
\left(I^\alpha_{t-}f\right)(x)=\frac{1}{\Gamma(\alpha)}\int_x^t f(s)(s-x)^{\alpha-1}ds.
\end{equation*}
The fractional derivatives are defined as follows.
\begin{maar}
Let $\alpha\in (0,1)$ and $f\in I_+^\alpha(L^1([0,t]))$. Then
\begin{equation*}
\left(D^\alpha_{0+}f\right)(x)=\frac{1}{\Gamma(1-\alpha)}\frac{d}{dx}\int_0^x f(s)(x-s)^{-\alpha}ds.
\end{equation*}
If $f\in I_-^\alpha(L^1([0,t]))$, then
\begin{equation*}
\left(D^\alpha_{t-}f\right)(x)=\frac{-1}{\Gamma(1-\alpha)}\frac{d}{dx}\int_x^tf(s)(s-x)^{-\alpha}ds.
\end{equation*}
\end{maar}
\begin{theor}
\label{integralexists}
Let $f\in W_2^\beta([0,T])$ and $g\in W_1^{1-\beta}([0,T])$. Then for any $t\in(0,T]$ there exists Lebesgue integral
\begin{equation*}
\int_0^t \left(D_{0+}^\beta f\right)(x)\left(D_{t-}^{1-\beta} g_{t-}\right)(x)dx.
\end{equation*}
\end{theor}
\begin{maar}[Generalized Lebesgue-Stieltjes integral]
We call
\begin{equation*}
\int_0^t fdg:=\int_0^t \left(D_{0+}^\beta f\right)(x)\left(D_{t-}^{1-\beta} g_{t-}\right)(x)dx
\end{equation*}
the generalized Lebesgue-Stieltjes integral of $f$ w.r.t. $g$.
\end{maar}
Note that the integral is the same for all $\beta$ for which it can be defined.
\begin{theor}
\label{tool}
Let $f,(f_n)_{n=1}^\infty\in W_2^\beta([0,T])$ and $g\in W_1^{1-\beta}([0,T])$. If
\begin{equation*}
\left|\left|f_n-f\right|\right|_{2,\beta}\rightarrow 0,
\end{equation*}
then
\begin{equation*}
\int_0^t f_ndg\rightarrow \int_0^t fdg,
\end{equation*}
for all $t\in(0,T]$.
\end{theor}
For more details, see~\cite{mishura}~or~\cite{a-m-v} and references therein.

\section*{Appendix B: Vertical and horizontal derivatives}
The following definitions and notations are taken from~\cite{cont}. We denote by $x_t$ the whole path of $x\in C([0,T])$ up to time $t$, that is $(x(u), 0\leq u\leq t)$. The vertical perturbation of path $x_t$ is defined for $h\in \R$ as
\begin{equation*}
x^h_t(u)= x(u), \quad u \in[0,t)
\end{equation*}
and
\begin{equation*}
x^h_t(t)=x(t)+h.
\end{equation*}
The horizontal extension of $x_t$ for $h>0$ is defined as
\begin{equation*}
x_{t,h}(u)=x(u), \quad u\in [0,t]
\end{equation*}
and
\begin{equation*}
x_{t,h}(u)=x(t), \quad u\in (t,t+h].
\end{equation*}
We say that a  family of maps $F=(F_t)_{t\in[0,T]}$, $F_t:C([0,t])\rightarrow \R$ is a non-anticipative functional. For the measurability issues we refer to~\cite{cont}. Now the horizontal derivative of $F$ at $x\in C([0,T])$ is defined as
\begin{equation}
\label{horizontal}
\mathcal{D}_tF(x)=\lim_{h\downarrow 0}\frac{F_{t+h}(x_{t,h})-F_t(x)}{h},
\end{equation} 
if the limit exists. If the limit of~(\ref{horizontal}) exists for all $x\in \bigcup_{t\in[0,T]}C([0,t])$, then the map
\begin{equation*}
\mathcal{D}_tF:C([0,t])\mapsto \R, \quad x\mapsto \mathcal{D}_t F(x)
\end{equation*}
defines a non-anticipative functional $\mathcal{D}F=(\mathcal{D}_tF)_{t\in[0,T)}$ that is called the horizontal derivative of $F$.

The vertical derivative is defined in the following way. A non-anticipative functional $F$ is vertically differentiable at $x\in C([0,t])$ if limit
\begin{equation}
\label{vertical}
\partial_xF_t(x)=\lim_{h\rightarrow 0}\frac{F_t(x^h_t)-F_t(x)}{h}
\end{equation}
exists. The limit of equation~(\ref{vertical}) is called the vertical derivative of $F$ at $x$. If the limit is defined for all $x\in \bigcup_{t\in[0,T]}C([0,t])$ then
\begin{equation*}
\partial_xF:C([0,t])\mapsto \R, \quad x\mapsto \partial_x F_t(x)
\end{equation*}
defines a non-anticipative functional $\partial_xF=(\partial_xF_t)_{t\in[0,T]}$.

We have the following chain rule for the vertical derivative
\begin{lemma}
\label{verticalchain}
Let $\phi \in C^1$ and $F$ be a non-anticipative functional. If $x\in C([0,t])$ is such that $\partial_x F_t(x)$ exists, then it holds for the vertical derivative that
\begin{equation*}
\partial_x \phi(F_t(x))=\phi'(F_t(x))\cdot \partial_x F_t(x).
\end{equation*}
\end{lemma}
\begin{proof}
The proof goes along the lines of the ordinary one dimensional chain rule. Note that
\begin{equation*}
F_t(x^h_t)-F_t(x)=h \partial_x F_t(x)+m(h)h,
\end{equation*}
where $m(h)\rightarrow 0$, when $h \rightarrow 0$. We also have that
\begin{equation*}
\phi(F_t(x)+\gamma)-\phi(F_t(x))=\gamma \phi'(F_t(x))+\alpha(\gamma)\gamma,
\end{equation*}
where $\alpha(\gamma)\rightarrow 0$ when $\gamma\rightarrow 0$. Now we have that
\begin{align*}
&\phi(F_t(x^h_t))-\phi(F_t(x))=\phi(F_t(x)+h \partial_x F_t(x)+m(h)h)-\phi(F_t(x))\\
=&(h\partial_xF_t(x)+m(h)h)\phi'(F_t(x))+\alpha(h\partial_xF_t(x)+m(h)h)(h\partial_xF_t(x)+m(h)h).
\end{align*}
If we now divide both sides by $h$ and let $h\rightarrow 0$ we obtain the claim.
\end{proof}
Next we will prove a product rule for the vertical derivative.
\begin{lemma}
\label{productrule}
Let $F$ and $E$ be non-anticipative functionals. Let $x\in C([0,t])$ such that $\partial_x F_t(x)$ and $\partial_x E_t(x)$ exist. Then it holds that
\begin{align*}
\partial_x \left(F_t\left(x\right)E_t\left(x\right)\right)
=F_t\left(x\right)\partial_x E_t\left(x\right)+E_t\left(x\right)\partial_x F_t\left(x\right).
\end{align*}
\end{lemma}
\begin{proof}
We have
\begin{align*}
&F_t\left(x^h_t\right)E_t\left(x^h_t\right)-F_t(x)E_t(x)\\
=&F_t\left(x^h_t\right)\left(E_t(x^h_t)-E_t(x)\right )+E_t(x)\left(F_t\left(x^h_t\right)-F_t(x)\right).
\end{align*}
Note that $\lim_{h\rightarrow 0}F_t(x_t^h)=F_t(x)$ and the vertical derivatives of $F$ and $E$ exist. Thus, the claim follows.
\end{proof}
For the horizontal derivative we have the following form of chain rule
\begin{lemma}
\label{horizontalchain}
Let $\phi \in C^1$ and $F$ be a non-anticipative functional. If $\mathcal{D}_tF(x)$ exists, then it holds for the horizontal derivative that
\begin{equation*}
\mathcal{D}_t \phi(F_t(x))=\phi'(F_t(x))\cdot \mathcal{D}_t F_t(x).
\end{equation*}
\end{lemma}
\begin{proof}
The proof of the lemma is analogous to the proof of lemma~\ref{verticalchain} and it is therefore omitted.
\end{proof}
\section*{Acknowledgements}
I have been supported financially by Academy of Finland, grant 21245. I am grateful to Ehsan Azmoodeh, Esko Valkeila and Lauri Viitasaari for their comments.

\end{document}